\documentclass[fleqn]{article}
 \usepackage{amscd}
                     \usepackage{amssymb}
 \newcommand{\be}{\begin{equation}}
       \newcommand{\ee}{\end{equation}}
       \newcommand{\ba}{\begin{eqnarray}}
        \newcommand{\ea}{\end{eqnarray}}
 \newcommand{\ban}{\begin{eqnarray*}}
 \newcommand{\ean}{\end{eqnarray*}}

 \newcommand{\lp}{\langle}
 \newcommand{\rp}{\rangle}
 \newcommand{\ra}{\rightarrow}

  \newcommand{\qed}{\hspace*{\fill}\rule{3mm}{3mm}\quad \vspace{.2cm}}
  \newcommand{\Pf}{\noindent {\bf Proof:} }
  \newcommand{\Rk}{\noindent {\bf Remark} }
 
 \renewcommand{\theequation}{\arabic{section}.\arabic{equation}}
 \newcommand{\sect}[1]{\section{#1} \setcounter{equation}{0}}

 \newcommand{\vol}{\mathrm{Vol}}
 \newcommand{\diam}{\mathrm{diam}}
 \newcommand{\Ric}{\mathrm{Ric}}
  \newcommand{\sn}{\mathrm{sn}}
 \newcommand{\Hess}{\mathrm{Hess}}

 \newtheorem{theo}{Theorem}[section]

\begin{document}
 \newtheorem{defn}[theo]{Definition}
 \newtheorem{ques}[theo]{Question}
 \newtheorem{lem}[theo]{Lemma}
 \newtheorem{lemma}[theo]{Lemma}
 \newtheorem{prop}[theo]{Proposition}
 \newtheorem{coro}[theo]{Corollary}
 \newtheorem{ex}[theo]{Example}
 \newtheorem{note}[theo]{Note}
 \newtheorem{conj}[theo]{Conjecture}
 \newtheorem{Ex}[theo]{Example}

\title{Comparison Geometry for the Bakry-Emery Ricci Tensor}
\author{Guofang Wei \thanks{Partially supported by NSF grant DMS-0505733}\and Will Wylie}
\date{}
\maketitle
\begin{abstract}
For Riemannian manifolds with a measure $(M,g, e^{-f} dvol_g)$ we
prove mean curvature and volume comparison results when the
$\infty$-Bakry-Emery Ricci tensor is bounded from below and $f$ is bounded or
$\partial_r f$ is bounded from below, generalizing the classical ones (i.e. when $f$
is constant).  This leads to extensions of many theorems for Ricci curvature bounded below to the Bakry-Emery Ricci tensor.  In particular, we give extensions of all of the major comparison theorems when $f$ is bounded.  Simple examples show the bound on $f$ is necessary for these results. 

\end{abstract}

\sect{Introduction} In this paper we study smooth metric measure
spaces $(M^n,g, e^{-f} dvol_g)$, where $M$ is a complete
$n$-dimensional Riemannian manifold with metric $g$, $f$ is a smooth
real valued function on $M$, and $dvol_g$ is the Riemannian volume
density on $M$. In this paper by the Bakry-Emery Ricci tensor we mean
$$\Ric_f = \Ric +\Hess f.$$
This is often also referred to as the $\infty$-Bakry-Emery Ricci
Tensor.    Bakry and Emery \cite{Bakry-Emery1985}  extensivley studied (and generalized) this tensor and its  relationship to diffusion processes.    The Bakry-Emery tensor  also occurs naturally in many
different subjects, see e.g. \cite{Lott2003} and
\cite[1.3]{Perelman-math.DG/0211159}. The equation $\Ric_f = \lambda
g$ for some constant $\lambda$ is exactly the gradient Ricci soliton
equation, which plays an important role in the theory of Ricci flow.
Moreover  $\Ric_f$ has a natural extension to metric measure spaces
\cite{Lott-Villani2007, Sturm2006-I, Sturm2006-II}. 

When $f$ is a
constant function, the Bakry-Emery Ricci tensor is  the Ricci tensor
so  it is  natural to investigate what geometric and topological
results for the Ricci tensor extend to the Bakry-Emery Ricci tensor.  Interestingly,  Lichnerowicz \cite{Lichnerowicz1970}  studied this problem at least 10 years before the work of Bakry and Emery. This  has also been actively investigated recently and there are a number
of interesting results in this direction which we will discuss
below. Also see  \cite{Chang-Gursky-Yang2006}
for another modification of the Ricci tensor and
Appendix~\ref{N-dim-appendix} for a discussion of  the
$N$-Bakry-Emery Ricci tensor $\Ric_f^N$ (see (\ref{Ric-N-f}) for the
definition). We thank Matthew Gursky for making us aware of \cite{Lichnerowicz1970}.

Although there is a Bochner formula for the Bakry-Emery Ricci
tensor \cite{Lott2003}  (see also (\ref{Bochner-f})) many of the
other basic geometric tools  for Ricci curvature do not extend. In
Section~\ref{examples} we give  a quick overview with examples
showing that, in general,  Myers' theorem, Bishop-Gromov's volume
comparison,  Cheeger-Gromoll's splitting theorem, and
Abresch-Gromoll's excess estimate are not true for the Bakry-Emery
Ricci tensor.  In this paper we find  conditions on $f$ that imply versions of these theorems.  In particular, we   show that all these theorems  hold when $f$ is
bounded.\footnote{ After writing the original version of this paper,  we learned that Lichnerowicz had proven the splitting theorem for $f$ bounded. We think this result is very interesting and  does not seem to be well known in the literature, so we have retained our complete proof here. }  These results give new tools for studying the Bakry-Emery
tensor and lead to generalizations of many of the classical
topological and geometric theorems for manifolds with a lower Ricci
curvature bound, and generalize all previous topological results for
the Bakry-Emery tensor.

For Ricci curvature all of the theorems listed above can be proven from the mean curvature (or Laplacian) comparison, see \cite{Zhu1997}.  Recall  that the mean curvature measures the relative rate of change of the volume element. Therefore, for the measure $e^{-f}dvol$, the associated mean curvature is $m_f =   m-\partial_r f,$  where $m$ is the mean curvature of the geodesic sphere with inward pointing normal vector.   Note that $m_f = \Delta_f (r)$, where $r$ is the distance function and  $\Delta_f = \Delta - \nabla f \cdot
 \nabla$ is the naturally associated ($f$-)Laplacian which is self-adjoint with respect to the weighted measure.  
 
In this paper we prove three mean curvature comparisons.  The first (see Theorem \ref{BasicMean}) is the most general and is quite simple to prove. Still, it has some interesting applications for manifolds with positive Bakry-Emery tensor (Corollaries \ref{finVol} and \ref{Liouville}). The other two are more delicate and have many applications.

\begin{theo}[Mean Curvature Comparison.] \label{mean-comp} Let $p \in M^n$ Assume $\Ric_f (\partial_r, \partial_r) \ge (n-1)H$,

a) if $\partial_r f \ge -a \ (a \ge 0)$ along a minimal geodesic
segment from $p$ (when $H>0$ assume $r \le \pi/2\sqrt{H}$) then
 \be \label{mean-a}
m_f(r) - m_H(r) \le a
 \ee
 along that minimal geodesic segment from $p$. Equality holds if and only if the radial sectional
 curvatures are equal to $H$ and $f(t) = f(p) - at$ for all $t<r$.

 b) if $|f| \le k$ along a minimal geodesic
segment from $p$ (when $H>0$ assume $r \le \pi/4\sqrt{H}$) then
 \be  \label{mean-H<0}
m_f(r)\le   m^{n+4k}_H(r)
  \ee
  along that minimal geodesic segment from $p$.  In particular when $H=0$ we have
\be \label{mean-H=0} m_f(r) \le \frac{n+4k-1}{r} \ee
\end{theo}  

Here $m^{n+4k}_H$ is the mean curvature of the
geodesic sphere in $M_H^{n+4k}$, the simply connected model space of dimension $n+4k$ with constant curvature $H$ and $m_H$ is the mean curvature of the model space of  dimension $n$.  See  (\ref{m-f-H>0}) in Section~\ref{mean-curvature} for the case $H>0$ and  $r\in [\frac{\pi}{4\sqrt{H}}, \frac{\pi}{2\sqrt{H}}]$ in part b. 

As in the classical case, these  mean curvature comparisons have  many applications.  First, we have volume comparison theorems.  Let $\vol_f (B(p,r)) = \int_{B(p,r)} e^{-f} dvol_g$, the weighted (or
 $f$-)volume and $\vol^n_H(r)$ be the volume of the radius $r$-ball in the model
space $M_H^n$. 

\begin{theo}[Volume Comparison.] \label{vol-comp} Let $(M^n,g, e^{-f}dvol_g)$ be  complete smooth metric measure space with $\Ric_f \geq
(n-1)H$. Fix $p \in M^n$.

a)  If $\partial_r f \ge -a$ along all
minimal geodesic segments from $p$ then for $R\ge r>0$ (assume $R \le
\pi/2\sqrt{H}$ if $H>0$) ,
 \be \label{vol-1} \frac{\vol_f (B(p,R))}{\vol_f (B(p,r))} \le e^{aR}
 \frac{\vol^n_H(R)}{\vol^n_H(r)}.
 \ee

  Moreover, equality holds if and only if the radial sectional
 curvatures are equal to $H$ and $\partial_rf \equiv a$.  In particular if $\partial_r f \geq 0$ and $\Ric_f \geq 0$ then $M$ has $f$-volume growth of degree at most $n$.

 b) If  $|f(x)| \leq k$ then for $R \ge r>0$ (assume $R \le
\pi/4\sqrt{H}$ if $H>0$), \be   \frac{\vol_f (B(p,R))}{\vol_f
(B(p,r))} \le
 \frac{\vol^{n+4k}_H(R)}{\vol^{n+4k}_H(r)}.
 \ee
 In particular, if $f$ is bounded and $\Ric_f \geq 0$ then $M$ has polynomial $f$-volume growth.
\end{theo}

\Rk \textbf{1} When $\Ric_f \geq 0$  the condition $f$ is bounded or
$\partial_r f \geq 0$ is necessary to show polynomial $f$-volume
growth as shown by Example~\ref{example-R}.  Similar statements are true for the volume of tubular
neighborhood of a hypersurface. See Section~{\ref{Volume} for
another version of volume comparison which holds for all $r>0$ even
when $H>0$.

\Rk \textbf{2} To prove the theorem we only need a lower bound on
$\Ric_f$  along the radial directions. Given any manifold $M^n$ with
Ricci curvature bounded from below one can always choose suitable
$f$ to get any lower bound for $\Ric_f$ along the radial directions.
For example if  $\Ric \ge -1$ and  $p \in M$, if we choose $f(x) =
r^2 = d^2(p,x)$, then $\Ric_f (\partial_r, \partial_r) \ge 1$.  Also
see Example~\ref{local-perturbation}.


\Rk \textbf{3} Volume comparison theorems have been proven for
manifolds with $N$-Bakry Emery Ricci tensor bounded below.  See Qian
\cite{Qian1997}, Bakry-Qian \cite{Bakry-Qian2005}, Lott
\cite{Lott2003}, and Appendix A. The $N$-Bakry Emery Ricci tensor is
\be \Ric_f^N = \Ric_f -\frac 1N df \otimes df  \ \ \  \mbox{for} \
N>0. \label{Ric-N-f} \ee For example, Qian shows that if $\Ric_f^N
\ge 0$ then $\vol_f(B(p,r))$ is of polynomial growth with degree
$\le n+N$. Note that $\Ric_f = \Ric_f^\infty$ so  one does not
expect polynomial volume growth for $\Ric_f \geq 0$.    Since
$\Ric_f^N \ge 0$ implies $\Ric_f \ge 0$ our result greatly improves
the volume comparison result of Qian when  $N$ is big and $f$ is
bounded, or when $\partial_r f \geq 0.$

The mean curvature and volume comparison theorems have many other applications.  We highlight two extensions of theorems of Calabi-Yau \cite{Yau1976} and Myers'  to the case where $f$ is bounded.    

\begin{theo} \label{inf-vol} If $M$ is a noncompact, complete manifold with $\Ric_f \geq 0$ for some bounded $f$
 then $M$ has at least linear $f$-volume growth.
\end{theo}

\begin{theo}[Myers' Theorem] \label{Myers'} If $M$ has $\Ric_f \geq (n-1)H > 0$ and $|f| \le k$ then $M$ is
compact and $\diam_M \leq \frac{\pi}{\sqrt{H}} + \frac{4k}{(n-1)
\sqrt{H}}$.
\end{theo}

 Examples~\ref{R^n} and \ref{example-H}
show that the assumption of bounded $f$ is necessary in both theorems. Qian \cite{Qian1997} has proven versions of both theorems for $\Ric_f^N$.  For other Myers' theorems see \cite{ Fern-Garcia, Zhang2007, Li-Xue-Mei1995, Morgan2006}.

The paper is organized as follows. In the next section we state and prove the mean curvature
comparisons.  In Sections~\ref{Volume} and \ref{Vol-appl} we prove the volume comparison theorems and discuss their applications, including Theorem \ref{inf-vol}.  In Section \ref{Splitting-Sec} we apply the mean curvature comparison to prove the splitting theorem for the Bakry-Emery tensor that is originally due to Lichnerowicz.  In Section \ref{Mean-appl} we discuss some other applications of the mean curvature comparison including the Myers' theorem and an  extension Abresch-Gromoll's excess estimate to $\Ric_f$. In Section 7 we discuss examples and questions.
Finally in Appendix A we state the mean curvature comparison for
$\Ric_f^N$.  This is a special case of an estimate in \cite{Bakry-Qian2005}, but we have written the result in more Riemannian geometry friendly language.  This gives other proofs of  the comparison theorems  for $\Ric_f^N$ mentioned above.

After posting the original version of this paper  we learned from
Fang, Li, and Zhang about their work which is closely related to
some of our work here. We thank them for sharing their work with us.
Their paper is now posted, see \cite{Fang-Li-Zhang}.  Motivated from
their paper we were able to strengthen the original version of
Theorem \ref{mean-comp} and Theorem \ref{vol-comp} and give a new
proof to Theorem~\ref{mean-comp}. This proof of the mean curvature
comparison seems to us to be new even in the classical Ricci
curvature case. We have moved our original proof using ODE methods
to an appendix because we feel it might be useful in other
applications.

 From the work of \cite{Petersen-Wei1997} one expects that the volume comparison and splitting theorem can be extended to the case that $\Ric_f$ is bounded from below in the integral sense. We also expect similar versions for metric measure spaces. These will be treated in  separate paper.

\textbf{Acknowledgment:} The authors would like to thank John Lott, Peter Petersen and  Burkhard Wilking for their interest and helpful discussions.

   \sect{Mean
Curvature Comparisons}  \label{mean-curvature}

 In this section we prove the mean curvature comparison theorems.   First we give a rough estimate on $m_f$ which is useful when $\Ric_f \ge \lambda g$
and $\lambda >0$.
\begin{theo} [Mean Curvature Comparison I.] \label{BasicMean} If $\Ric_f (\partial_r, \partial_r) \ge \lambda$ then given any minimal geodesic segment and $r_0 >0$,
 \be \label{m_H>0}
m_f(r) \le m_f(r_0) - \lambda (r-r_0) \ \  \mbox{for}\ r \ge r_0.
 \ee
Equality holds for some $r>r_0$ if and only if  all the radial sectional curvatures are zero, $\Hess r\equiv 0$, and $\partial_r^2 f \equiv \lambda$ along the geodesic from $r_0$ to $r$.
 \end{theo}

 \Pf Applying  the Bochner formula
 \be \label{bochner} \frac 12 \Delta
|\nabla u|^2 = |\Hess \, u|^2 + \lp \nabla u, \nabla (\Delta u) \rp
+ \Ric (\nabla u, \nabla u)
 \ee to the
distance function $r(x) = d(x,p)$, we have
\be \label{bochner-distance}
0 = | \Hess \, r |^2 + \frac{\partial}{\partial r} (\Delta r) + \Ric
(\nabla r, \nabla r).
\ee
Since $\Hess \, r$ is the second fundamental from of the geodesic
sphere and $\Delta r$ is the mean curvature, with the Schwarz
inequality, we have the Riccati inequality
 \be \label{Riccati-inequ} m' \le -\frac{m^2}{n-1} -
\Ric (\partial r, \partial r). \ee And equality holds if and only if
the radial sectional
 curvatures are constant.
Since $m_f' = m'- \Hess f \, (\partial r, \partial r)$, we have \be
\label{Riccati-inequ-mf} m_f' \le -\frac{m^2}{n-1} - \Ric_f
(\partial r, \partial r). \ee If $\Ric_f \ge \lambda$, we have
\[
m_f' \le -\lambda.\] This immediately gives the inequality
(\ref{m_H>0}).

To see the equality statement, suppose $m_f ' \equiv -\lambda$ on an
interval $[r_0, r]$, then from (\ref{Riccati-inequ-mf}) we have $m
\equiv 0$  and \be m_f' = -\partial_r^2 f = -\Ric_f(\partial_r ,
\partial_r )=-\lambda.\ee So we also have $\Ric( \partial_r,
\partial_r) = 0$.  Then by (\ref{bochner-distance})  $\Hess \, r =
0$, which implies the sectional curvatures must be zero. \qed
\newline

\noindent {\bf Proof of Theorem~\ref{mean-comp}}.  Let $\sn_H(r)$ be the solution to \[ \sn_H''+H \sn_H = 0\] such
that $\sn_H(0)=0$ and $\sn_H'(0)=1$. Then \be \label{relation}
m^n_H = (n-1) \frac{\sn'_H}{\sn_H}. \ee So we have \ba
\left(\sn^2_Hm\right)'  & =& 2 \sn'_H \sn_H m + sn^2_H m' \\
& \leq& 2 \sn'_H \sn_H m - \frac{\sn^2_H m^2}{n-1} - \sn^2_H \Ric
(\partial_r, \partial_r)\\ & = & -\left(\frac{\sn_H m}{\sqrt{n-1}} -
\sqrt{n-1} \sn'_H \right)^2 + (n-1)(\sn'_H)^2  -
\sn^2_H\Ric(\partial_r, \partial_r)   \\ & \leq & (n-1) (\sn'_H)^2 -
(n-1) H\sn^2_H + \sn^2_H
\partial_r \partial_r f. \ea Here in the 2nd line we have used
(\ref{Riccati-inequ}), and in the last we used the lower bound on
$\Ric_f$.

On the other hand (\ref{relation}) implies that \[ (\sn^2_H m_H)' =
(n-1)(\sn'_H)^2 - (n-1)H\sn^2_H. \] Therefore we have \be
\left(\sn^2_Hm\right)' \leq \left(\sn^2_Hm_H\right)' + \sn^2_H
\partial_t \partial_t f . \ee
 Integrating from 0 to r yields \be
\sn^2_H(r) m(r) \leq \sn^2_H(r) m_H(r) + \int_{0}^r \sn^2_H(t)
\partial_t \partial_t f(t) dt.\ee
When $f$ is constant (the classical case) this gives the usual mean
curvature comparison. This quick proof does not seem to be in the
literature.

 \textbf{Proof or Part a.} Using integration by parts on the last
term we have \be \label{PartsOnce} \sn^2_H(r) m_f(r) \leq \sn^2_H(r)
m_H(r) - \int_{0}^r (\sn^2_H(t))'
\partial_t f(t) dt.\ee Under our assumptions $(\sn^2_H(t))' = 2
\sn'_H(t) \sn_H(t)\geq 0$  so if $\partial_t f(t) \geq -a$ we have \be
\sn^2_H(r) m_f(r) \leq \sn^2_H(r) m_H(r) + a\int_{0}^r (\sn^2_H(t))'
dt = \sn^2_H(r) m_H(r)+ \sn^2_H(r) a\ee This
proves the inequality.

To see the rigidity statement suppose that $\partial_t f \geq -a$
and $m_f(r) = m_H(r) + a$ for some $r$.  Then from (\ref{PartsOnce})
we see \be a \sn^2_H \leq \int_{0}^r (\sn^2_H(t))'
\partial_t f(t) dt \leq a \sn^2_H. \ee So that $\partial_t f \equiv
-a$.  But then $m(r) = m_f - a = m_H(r)$ so that the rigidity
follows from the rigidity for the usual mean curvature comparison.

\textbf{Proof of Part b.} Integrate  (\ref{PartsOnce}) by parts
again \be  \label{Better0}\sn^2_H(r) m_f(r) \leq \sn^2_H(r) m_H(r) -
f(r) (\sn^2_H(r))' +  \int_{0}^r f(t) (sn^2_H)''(t) dt. \ee Now if
$|f| \leq k$ and $r \in (0, \frac{\pi}{4\sqrt{H}}]$ when $H>0$, then
$(sn^2_H)''(t) \ge 0$ and we have \be \sn^2_H(r) m_f(r) \leq
\sn^2_H(r) m_H(r) + 2k (\sn^2_H(r))'. \ee From (\ref{relation}) we
can see that
\[ (\sn^2_H(r))' = 2 \sn'_H\sn_H = \frac{2}{n-1}m_H \sn^2_H. \] so we
have \be m_f(r) \leq \left(1+ \frac{4k}{n-1}\right) m_H(r) =
m^{n+4k}_H(r). \ee \qed

 Now when $H>0$ and $r \in
[\frac{\pi}{4\sqrt{H}}, \frac{\pi}{2\sqrt{H}}]$, \ban
 \int_{0}^r f(t) (sn^2_H)''(t) dt & \le & k \left( \int_{0}^{\frac{\pi}{4\sqrt{H}}}(sn^2_H)''(t) dt
 - \int_{\frac{\pi}{4\sqrt{H}}}^{r}(sn^2_H)''(t) dt\right) \\
 & = & k \left( \frac{2}{\sqrt{H}} - sn_H(2r) \right).
 \ean
 Hence
 \be \label{m-f-H>0}
  m_f(r) \leq \left(1+ \frac{4k}{n-1}\cdot \frac{1}{\sin (2\sqrt{H} r)} \right)
m_H(r). \ee This estimate will be used later to prove the Myers'
theorem in Section 5.

 \Rk In the case $H=0$, we have $sn_H (r) = r$ so
(\ref{Better0}) gives the estimate in \cite{Fang-Li-Zhang} that \be
\label{FLZ} m_f(r) \leq \frac{n-1}{r} - \frac{2}{r}f(r) +
\frac{2}{r^2} \int_{0}^r f(t) dt. \ee

\Rk  The exact same argument gives mean curvature comparison for the mean curvature of distance sphere of hypersurfaces with $\Ric_f$ lower bound.

\sect{Volume Comparisons} \label{Volume}
In this section we prove
the  volume comparison theorems.

For $p \in M^n$, use exponential polar coordinate around $p$ and
write the volume element $d\, vol = \mathcal A(r,\theta) dr \wedge
d\theta_{n-1}$, where $d\theta_{n-1}$ is the standard volume element
on the unit sphere $S^{n-1}(1)$. Let $\mathcal A_f(r, \theta) =
e^{-f}\mathcal A(r,\theta)$.  By the first variation of the area
(see \cite{Zhu1997})
\be \label{area} \frac{\mathcal A'}{\mathcal A}(r,\theta) = (\ln (\mathcal A(r,\theta)))' = m (r, \theta).
 \ee
Therefore \be \label{area-f} \frac{\mathcal A_f'}{\mathcal A_f}
(r,\theta) = (\ln (\mathcal A_f(r,\theta)))' = m_f (r, \theta). \ee
And for $r\ge r_0 >0$
\be \label{area-m}
\frac{\mathcal A_f(r,\theta)}{\mathcal A_f(r_0,\theta)} =
e^{\int_{r_0}^r m_f (r, \theta)}.\ee
The volume comparison theorems follow  from the mean curvature comparisons through this equation.

 First applying the mean curvature estimate Theorem~\ref{BasicMean} we have the following basic volume comparison theorem.

\begin{theo} [Volume Comparison I]  \label{BasicVol} Let $\Ric_f \geq \lambda$ then for any $r$ there are constants $A$, $B$, and $C$ such that \[ \vol_f(B(p,R)) \leq A + B \int_r^R e ^{-\lambda t^2 + Ct}dt.\]
\end{theo}

The  version of Theorem~\ref{BasicVol} for tubular neighborhoods of
hypersurfaces is very similar and has been  proven by Morgan
\cite{Morgan2005}. As Morgan points out, the  theorem is optimal and
the constants can not be uniform as the Gaussian soliton shows, see
Example~\ref{R^n}.

\Pf Using the mean curvature estimate (\ref{m_H>0})
\[
\int_{r_0}^r m_f(r) \le m_f(r_0)r - \frac 12 \lambda r^2.
\]
Hence
\[
\mathcal A_f(r,\theta) \le \mathcal A_f(r_0,\theta) e^{m_f(r_0,
\theta)r - \frac 12 \lambda r^2}.
\]
Now let $A(p,r_0,r)$ be the  the annulus $A(p,r_0,r) =
B(p,r)\setminus B(p,r_0)$.  Then \ba
\vol_f (A(p,r_0,r)) & = &  \int_{r_0}^r \int_{S^{n-1}} \mathcal A_f(s,\theta) d\theta ds  \\
& \le & \int_{r_0}^r \int_{S^{n-1}} \mathcal A_f(r_0,\theta) e^{m_f(r_0,\theta)r - \frac 12 \lambda r^2}d\theta \label{positive-vol-comp}ds \\
& \le & A_f(r_0) \int_{r_0}^r e^{Cr -\frac 12 (n-1)\lambda
r^2}ds.\ea Where $A_f(r_0)$ is the surface area of the geodesic
sphere induced from the $f$-volume element and $C$ is a
constant such that $C \geq m_f(r_0, \theta)$ for all $\theta$ where
it is defined. Since $\vol_f(B(p,r)) = \vol_f(\vol(B(p,r_0)) +
\vol_f (A(p,r_0,r))$ this proves the theorem. \qed

We also have a rigidity statement for (\ref{positive-vol-comp}).  That
is, if the inequality (\ref{positive-vol-comp}) is an equality then
we must have equalities in the mean curvature comparison along all the
geodesics, this implies that $\Hess \, r \equiv 0$ which implies
that \be A(p, r_0,r) \cong S(p,r_0) \times [r_0,r] \ee where
$S(p,r_0$) is the geodesic sphere with radius $r_0$.  Moreover,
$f(x,t)= f(x)+ \partial_r f(x)(r-r_0)+ \frac{\lambda}{2}(r-r_0)^2$.

Now we prove Theorem~\ref{vol-comp} using Theorem~\ref{mean-comp}.

\noindent {\bf Proof of Theorem \ref{vol-comp}}: For Part a) we
compare with a model space, however, we modify the measure according
to a.  Namely, the  model space will be the pointed metric measure
space $M^n_{H,a}=(M^n_H,g_H, e^{-h}dvol, O)$ where $(M^n_H,g_H)$ is
the n-dimensional simply connected space with constant sectional
curvature $H$, $O \in M^n_H$, and $h(x) = -a \cdot d(x,O)$.  We make
the model a pointed space because the space only has
$\Ric_f(\partial_r, \partial_r) \geq (n-1)H$ in the radial
directions from $O$ and we only compare volumes of balls centered at
$O$.

Let  ${\mathcal A}^a_H$ be the $h$-volume element in
$M_{H,a}^n$. Then ${\mathcal A}^a_H(r) = e^{ar} \mathcal
A_H(r)$ where $\mathcal
A_H$ is the Riemannian volume element in $M^n_H$. By the mean curvature comparison we have $(\ln(\mathcal{A}_f(r, \theta))'\leq a+m_H = (\ln ({\mathcal
A}_H^a))'$ so  for $r<R$, \be \label{vol-elem}
\frac{\mathcal A_f(R,\theta)}{\mathcal A_f(r,\theta)} \le
\frac{\mathcal A_H^{a}(R, \theta)}{\mathcal A_H^{a}(r,\theta)}.
\ee
Namely $\frac{\mathcal A_f(r,\theta)}{\mathcal
A_H^{a}(r,\theta)}$ is nonincreasing in $r$. Using Lemma~3.2 in
\cite{Zhu1997}, we get for $0< r_1 < r$, $0< R_1 <R$, $r_1 \le R_1$, $r
\le R$, \be \label{vol} \frac{\int_{R_1}^R \mathcal
A_f(t,\theta)dt}{\int_{r_1}^r \mathcal A_f(t,\theta)dt} \le
\frac{\int_{R_1}^R \mathcal A_H^{a}(t, \theta)dt}{\int_{r_1}^r
\mathcal A_H^{a}(t,\theta)dt} . \ee
Integrating  along the sphere direction
gives \be \label{ann-comp-b} \frac{\vol_f (A(p,R_1,R))}{\vol_f (A(p,r_1,r))} \le
 \frac{\vol^a_H(R_1,R)}{\vol^a_H(r_1,r)}.
 \ee
Where $\vol_H^a(r_1,r)$ is the $h$-volume of the annulus $B(O,r) \setminus B(O,r_1) \subset M_H^n$.
Since $\vol_H(r_1, r) \leq \vol_H^a(r_1,r) \leq e^{ar} \vol_H(r_1,r)$ this gives (\ref{vol-1}) when $r_1=R_1 = 0$ and proves Part b).

In the model space the radial function $h$ is not smooth at the origin.  However, clearly one can  smooth the function to a  function with $\partial_r h \geq -a$ and $\partial_r^2 h \geq 0$ such that the $h$-volume taken with the smoothed $h$ is arbitrary close to that of the model. Therefore, the inequality (\ref{ann-comp-b}) is optimal. Moreover, one can see from the equality case of the mean curvature comparison that if the annular volume is equal to the volume in the model then all the radial sectional curvatures are $H$ and $f$ is exactly a linear function.
\newline

\noindent {\bf Proof of Part b)}: In this case let $\mathcal{A}^{n+4k}_H$ be the volume element in the simply connected model space with constant curvature $H$ and dimension $n+4k$.

Then from the mean curvature comparison we have $\ln(\mathcal{A}_f (r, \theta))' \leq \ln(\mathcal{A}^{n+4k}_H(r))'$.
So again applying Lemma~3.2 in \cite{Zhu1997}  we
obtain \be \label{ann-vol-c}
\frac{\vol_f(A(p,R_1,R))}{\vol_f(A(p,r_1,r))} \leq \frac{Vol^{n+4k}_H(R_1,R)}{Vol^{n+4k}_H(r_1, r)}.
\ee With $r_1=R_1=0$ this implies the relative volume
comparison for balls \be \frac{\vol_f(B(p,R))}{\vol_f(B(p,r))} \leq
 \frac{Vol^{n+4k}_H(R)}{Vol^{n+4k}_H(r)}.   \label{vol-H=0} \ee
Equivalently
 \be \frac{\vol_f(B(p,R))}{V^{n+4k}_H(R)} \leq
\frac{\vol_f(B(p,r))}{V^{n+4k}_H(r)}. \ee Since $n+4k > n$
we note that the right hand side blows up as $r \to 0$ so one does not obtain a uniform upper bound on $\vol_f(B(p,R))$. Indeed, it is not possible to do so since one can always add a constant to $f$ and not effect the Bakry-Emery tensor.

By taking $r=1$ we do obtain a volume growth estimate for $R > 1$
\be \vol_f(B(p,R)) \leq \vol_f(B(p,1)) Vol^{n+4k}_H(R) . \ee  Note
that, from Part a)  $\vol_f(B(p,1)) \leq e^{-f(p)} e^a \omega_n$
if $\partial_r f \geq -a$ on $B(p,1)$. \qed

In the next section we collect the applications of the volume comparison theorems.

\sect{Applications of the volume comparison theorems.} \label{Vol-appl}
In the case where $\lambda > 0$  Theorem~\ref{BasicVol} gives two
very interesting corollaries.  The first is also observed in \cite{Morgan2005}.

\begin{coro}\label{finVol}   \label{fingroup}If M is complete and $Ric_f \ge \lambda > 0$ then $\vol_f(M)$ is finite and $M$ has finite fundamental group.
\end{coro}

We note the finiteness of volume is true in  the setting of more
general diffusion operators \cite{Bakry-Emery1985} and more general metric measure spaces \cite{Hinde}. Using a
different approach the second author has proven that the fundamental
group is  finite for spaces satisfying  $\Ric +
\mathcal{L}_Xg \geq \lambda >0$ for some vector field $X$ \cite{Wylie}.   This
had earlier been shown under the additional assumption that the
Ricci curvature is bounded by Zhang \cite{Zhang2007}. See also
\cite{Naber}.     When $M$ is
compact the finiteness of fundamental group was first shown by X. Li
\cite[Corollary  3]{Li-Xue-Mei1995} using a probabilistic method.
Also see \cite{ Zhang2007, Fern-Garcia, Qian1997, Lott2003}.  We
would like to thank Prof. David Elworthy for bringing the article
\cite{Li-Xue-Mei1995} to our attention.
 
 The second corollary is the following Liouville Theorem, which is a strengthening of a result of Naber \cite{Naber}.

\begin{coro} \label{Liouville} If M is complete with $Ric_f \ge \lambda > 0$, $u \geq 0$,
$\Delta_f(u) \geq 0$, and there is $\alpha < \lambda$ such that  $u(x) \leq e^{\alpha d(p,x)^2}$ for some
$p \in M$ then $u$ is constant. \end{coro}

In particular there are no bounded f-subharmonic functions.  Corollary \ref{Liouville} follows from Yau's proof that a complete manifold has no  positive $L^p$
($p>1$) subharmonic functions \cite{Yau1976}. The argument only uses integration
by parts and a clever choice of test function and so is valid also
for the weighted measure and Laplacian.  

While Theorem \ref{BasicVol} has applications when $\lambda >0$ it is not strong enough to extend results for a general lower bound, for these results we apply Theorem \ref{vol-comp}. It is well known that a lower bound on volume growth for manifolds with $\Ric \geq 0$ can be derived from the volume comparison for annulli, see \cite{Zhu1997}.  Thus Theorem \ref{inf-vol} follows from (\ref{ann-vol-c}).  We give the proof here for completeness and to motivate Theorem \ref{inf-vol-conv}.

\noindent \textbf{Proof of Theorem~\ref{inf-vol}}: Let $M$ be a
manifold with $\Ric_f \geq 0$ for a bounded function $f$. Let $p \in
M$ and let $\gamma$ be a geodesic ray based at $p$  in $M$.   Then,
applying the annulus relative volume comparison (\ref{ann-vol-c}) to
annuli centered at $\gamma(t)$, we obtain \[
\frac{\vol_f(B(\gamma(t), t-1))}{\vol_f(A(\gamma(t), t-1, t+1))}
\geq \frac{(t-1)^{n+4k}}{ (t+1)^{n+4k} - (t-1)^{n+4k }} \geq c(n,k)t
\qquad \forall \ t \ge 2. \]

But $B(\gamma(0), 1) \subset A(\gamma(t), t-1, t+1)$ so we have \[
\vol_f(B(p, t - 1)) \geq c(n,k) \vol_f(B(p,1)) t  \qquad \forall t
\ge 2. \]  \qed

Using the volume comparison (\ref{ann-comp-b}) in place of  (\ref{ann-vol-c}) we can also prove a lower bound on the volume growth for certain convex $f$.

\begin{theo} \label{inf-vol-conv} If $\Ric_f \geq 0$ where $f$ is convex function such that the set of critical points of $f$  is unbounded, then $M$ has at least linear $f$-volume growth.
\end{theo}

The hypothesis on the critical point  set is necessary by Examples \ref{R^n} and \ref{example-H}.

\Pf
Fix $p \in M$.  Since the set of critical points of a convex function is connected, for every $t$ there is $x(t)$, a critical point of $f$, such that $d(p,x(t)) = t$.   But  $\nabla f(x(t)) =0$ and $f$ is convex so $\partial r f \geq 0$ in all the radial directions from $x(t)$, therefore we can apply (\ref{ann-comp-b}) and repeat the arguments in the proof of Theorem ~\ref{inf-vol} to prove the result.  \qed

In \cite{Milnor1968} Milnor observed that  polynomial volume growth on the universal cover of a manifold restricts the structure of its fundamental group.  Thus Theorem \ref{vol-comp} also implies the following extension of Milnor's Theorem.

\begin{theo} \label{poly} Let  $M$ be a complete manifold with $Ric_f \geq 0$.
\begin{enumerate}
\item  If  $f$ is a convex function that obtains its minimum then any finitely generated subgroup of $\pi_1(M)$ has polynomial growth of degree less than or equal $n$.  In particular, $b_1(M) \leq n.$
\item  If $|f| \leq k$ then any finitely generated subgroup of $\pi_1(M)$ has polynomial growth of degree less than
 or equal to $n+4k$.  In particular, $b_1(M) \leq n+4k$.
\end{enumerate}
\end{theo}

Part 1 follows because at a pre-image of the minimum point in the universal cover, $\partial_r f \geq 0$. Gromov \cite{Gromov1981} has shown that a finitely generated group
has polynomial growth if and only if it is almost nilpotent.
Moreover, the work of the first author and Wilking shows that any
finitely generated almost nilpotent group is the fundamental group
of a manifold with $\Ric \geq 0$ \cite{Wei1988, Wilking2000}.
Therefore, there is a complete classification of  the finitely
generated groups that can be realized as the fundamental group of a
complete manifold with $\Ric \geq 0$.  Combining these results with
Theorem~\ref{poly} we expand this classification to a larger class
of manifolds.

\begin{coro} \label{nilpotent} A finitely generated group $G$ is the fundamental group of some manifold with
\begin{enumerate}
\item $\Ric _f \geq 0$ for some bounded $f$ or
\item $\Ric_f \geq 0$ for some convex $f$ which obtains its minimum
\end{enumerate}
if and only if $G$ is almost nilpotent.
\end{coro}

It would be interesting to know if  Corollary \ref{nilpotent} holds without any
assumption on $f$.  Example~\ref{example-R} shows that
the Milnor argument can not be applied since the $f$-volume growth of a manifold with $\Ric_f \geq0$
may be exponential, so a different method of proof would be needed.

In \cite{Anderson1990-2} Anderson uses similar covering arguments to show, for example, that if $\Ric \geq 0$ and and $M$ has euclidean volume growth then $\pi_1(M)$ is finite.  He also finds interesting relationships between the first betti number, volume growth, and finite generation of fundamental group of manifolds with $\Ric \geq 0$.  These relationships also carry over to manifolds satisfying the hypotheses of Theorem \ref{poly}.  We leave these statements to the interested reader.

Applying the relative volume comparison Theorem~\ref{vol-comp} we also 
have the following extensions of theorems of Gromov
\cite{Gromov1999} and Anderson \cite{Anderson1990}.
\begin{theo}
For the class of manifolds $M^n$ with $\Ric_f \geq (n-1)H,$ $\diam_M
\leq D$ and $|f| \le k \ (|\nabla f|  \le  a)$, the first Betti
number $b_1 \le C(n,k,HD^2)\ (C(n, HD^2, aD))$.
\end{theo}

\begin{theo}
For the class of manifolds $M^n$ with $\Ric_f \geq (n-1)H,$ $\vol_f
\geq V$, $\diam_M \leq D$ and $|f| \le k \ (|\nabla f|  \le a)$
there are only finitely many isomorphism types of $\pi_1 (M)$.
\end{theo}

\Rk In the case when $|\nabla f|$ is bounded, $\Ric_f$  bounded from
below implies $\Ric_f^N$ is also bounded from below (with different
lower bound). Therefore in this case the results can also been
proven using the volume comparison in \cite{Qian1997, Lott2003, Bakry-Qian2005} for the
$\Ric_f^N$ tensor.

\sect{The Splitting Theorem.} \label{Splitting-Sec}

An important application of the mean curvature comparison is the extension of  the Cheeger-Gromoll splitting theorem.  After writing the original version of this paper,  we learned that Lichnerowicz had proven the splitting theorem \cite{Lichnerowicz1970}.  

\begin{theo}[Lichnerowicz-Cheeger-Gromoll Splitting Theorem] \label{splitting}
If $\Ric_f \ge 0$ for some bounded $f$ and $M$ contains a line, then $M=
N^{n-1} \times \mathbb R$ and $f$ is constant along the line.
\end{theo}

Since  Lichnerowicz  did not write out a detailed proof, 
we retain our complete proof here. 

\Rk In \cite{Fang-Li-Zhang} Fang, Li, and Zhang  show that only an
upper bound on $f$ is needed in the above theorem.
Example~\ref{example-H} shows that the upper bound on $f$ is necessary.

 Recall that $m_f = \Delta_f (r)$, the $f$-Laplacian of the distance
function. From (\ref{mean-comp}), we get a local Laplacian
comparison for distance functions \be  \label{lap-com} \Delta_f (r)
\le \frac{n+4k-1}{r}   \ \ \mbox{for all}\ x \in M\setminus \{p,
C_p\}. \ee

Where $C_p$ is the cut locus of $p$. To prove the splitting theorem we apply this estimate to the Busemann functions.

\begin{defn} If $\gamma$ is a ray then Busemann function associated to $\gamma$ is the function
\be
 b^{\gamma}(x) = \lim_{t \rightarrow \infty} (t - d(x, \gamma(t))).  \ee
 \end{defn}

>From the triangle inequality the Busemann function is Lipschitz continuous with Lipschitz constant 1 and thus is differential almost everywhere. At the points where $b_{\gamma}$ is not smooth we interpret the $f$-laplacian in the sense of barriers.

\begin{defn} For a continuous function $h$ on $M, q \in M$, a function $h_q$
defined in a neighborhood $U$ of $q$, is a lower barrier of $h$ at
$q$ if $h_{q}$ is $C^2(U)$    and
 \be  \label{up-barrier}
h_{q} (q) = h(q), \ \ \ h_{q} (x) \le h(x)  \  (x \in U). \ee
\end{defn}

\begin{defn} We say that $\Delta_f(h) \geq a$ in the barrier sense if, for every $\varepsilon > 0$, there exists a lower barrier function $h_{\varepsilon}$ such that $\Delta_f(h_{\varepsilon}) > a - \varepsilon$. An upper bound on $\Delta_f$ is defined similarly in terms of upper barriers.
\end{defn}

The local Laplacian comparison is applied to give the following key
lemma.

\begin{lemma} \label{bus-lap} If $M$ is a complete, noncompact manifold with $\Ric_f \geq 0$ for some bounded function $f$ then $\Delta_f(b^{\gamma}) \geq 0$ in the barrier sense.
\end{lemma}

\Rk As in \cite{Fang-Li-Zhang}, one can use the inequality
(\ref{FLZ}) to prove  Lemma \ref{bus-lap} only assuming an upper
bound on  $f$.

\Pf For the Busemann function at a point $q$ we have a family of
barrier functions defined as follows. Given $t_i \rightarrow
\infty$,  let $\sigma_i$ be minimal geodesics from $q$ to
$\gamma(t_i)$, parametrized by arc length.  The sequence
$\sigma_i'(0)$ subconverges to some $v_0 \in T_q M$.  We call  the
geodesic $\overline{\gamma}$ such that $\overline{\gamma}'(0) = v_0$
an asymptotic ray to $\gamma$.

Define the function  $h_t(x) = t - d(x, \overline{\gamma}(t)) + b^{\gamma}(q)$. Since $\overline{\gamma}$ is a ray, the points $q = \overline{\gamma}(0)$ and $\overline{\gamma}(t)$ are not cut points to each other, therefore the function $d(x, \overline{\gamma}(t))$ is smooth in a neighborhood of $q$ and thus so is $h_t$.  Clearly $h_t(q) = b^{\gamma}(q)$, thus to show that $h_t$ is a lower barrier for $b^{\gamma}$ we only need to show that  $h_t(x) \leq   b^{\gamma}(x)$.  To see this, first note that for any $s$,
\be  - d(x, \overline{\gamma}(t)) \leq    - d(x, \gamma(s)) + d(\gamma(s), \overline{\gamma}(t)) = s - d(x, \gamma(s)) - s +  d(\gamma(s), \overline{\gamma}(t)). \ee
Taking $s \rightarrow \infty$ this gives
\be \label{sp1}  - d(x, \overline{\gamma}(t)) \leq  b^{\gamma}(x) - b^{\gamma}(\overline{\gamma}(t)). \ee
Also,
\ba b^{\gamma}(q) & = & \lim_{i \rightarrow \infty} (t_i - d(q, \gamma(t_i))) \nonumber \\ &  =  &  \lim_{i \rightarrow \infty} (t_i - d(q, \sigma_i(t)) - d(\sigma_i(t),\gamma(t_i))) \nonumber \\ &=& -d(q, \overline{\gamma}(t)) + \lim_{i \rightarrow \infty} (t_i - d(\sigma_i(t),\gamma(t_i)))  \nonumber \\ &=& \label{sp2}-t + b^{\gamma}(\overline{\gamma}(t)). \ea
Combining (\ref{sp1}) and (\ref{sp2}) gives
\be  \label{Barrier} h_t(x) \leq b^{\gamma}(x), \ee
so  $h_t$ is a lower barrier function for $b^{\gamma}$. By (\ref{lap-com}), we have that
\be \Delta_f(h_t)(x) = \Delta_f(- d(x, \overline{\gamma}(t))) \geq -\frac{n+4k-1}{t}.  \ee
Taking $t \rightarrow \infty$ proves the lemma.
\qed

Note that since $\Delta_f$ is just a perturbation of $\Delta$ by a lower order term, the strong maximum principle and elliptic regularity still hold for $\Delta_f$.  Namely if $h$ is a continuous function such that $\Delta_f(h) \geq 0$ in the barrier sense and $h$ has an interior maximum then $h$ is constant and if $\Delta_f(h) =0$ (i.e $\ge 0$ and $\le 0$) in the barrier sense then $h$ is smooth.  We now apply the lemma and these two theorems to finish the proof of the splitting theorem.

\noindent {\bf Proof of Theorem \ref{splitting}}:
Denote by $\gamma_+$ and $\gamma_-$ the two rays which form the line $\gamma$ and let $b^{+}$, $b^{-}$ denote their Busemann functions.

The function $b^{+} + b^{-}$ has a maximum at $\gamma(0)$ and satisfies $\Delta_f (b^{+} + b^{-}) \geq 0$, thus by the strong maximum principle the function is constant and equal to 0.  But then $b^{+} = - b^{-}$ so that $0 \leq \Delta_f(b^{+})=-\Delta_f(b^{-}) \leq  0$ which then implies, by elliptic regularity, that the functions $b^{+}$ and $b^{-}$ are smooth.

 Moreover, for any point $q$ we can consider asymptotic rays $\overline{\gamma}_+$ and $\overline{\gamma}_-$ to $\gamma_+$and $\gamma_-$ and denote their Busemann functions by $\overline{b}^+$ and $\overline{b}^-$.  From (\ref{Barrier}) it follows that
 \be \overline{b}^+(x) + b^+(q) \leq b^+(x). \label{half} \ee  We will show that this inequality is, in fact, an equality when $\gamma_+$ extends to a line.

 First we show that the two asymptotic rays $\overline{\gamma}_+$ and $\overline{\gamma}_-$ form a line.  By the triangle inequality, for any t
 \ban
 d(\overline{\gamma}_+(s_1), \overline{\gamma}_-(s_2)) &\geq& d( \overline{\gamma}_-(s_2), \gamma_+(t)) - d(\gamma_+(t), \overline{\gamma}_+(s_1)) \\ & =& t - d(\gamma_+(t), \overline{\gamma}_+(s_1)) -(t -  d( \overline{\gamma}_-(s_2), \gamma_+(t)),  \ean
 so by taking $t \rightarrow \infty$ we have
\ban d(\overline{\gamma}_+(s_1), \overline{\gamma}_-(s_2)) & \geq & b^+(\overline{\gamma}^+(s_1) ) -b^+(\overline{\gamma}^-(s_2)) \\ & = &  b^+(\overline{\gamma}^+(s_1) ) +b^-(\overline{\gamma}^-(s_2)) \\ &\geq& \overline{b}^+(\overline{\gamma}^+(s_1) ) + b^+(q) + \overline{b}^-(\overline{\gamma}^-(s_2))+b^-(q)\\ & = & s_1+s_2. \ean

Thus, any asymptotic ray to $\gamma_+$ forms a line with any asymptotic ray to $\gamma_{-}$.    Applying the same argument given above for $b^+$ and $b^-$  we see that $\overline{b}^+ = - \overline{b}^-$.  By Applying  (\ref{half}) to $b^-$   \[-\overline{b}^-(x) - b^-(q) \geq -b^-(x). \]
 Substituting $b^+ = - b^-$ and $\overline{b}^+ = - \overline{b}^-$ we have
\[ \overline{b}^+(x) +b^+(q) \geq b^+(x). \] This along with (\ref{half}), gives
\[ \overline{b}^+(x) +b^+(q) = b^+(x). \]
Thus, $\overline{b}^+$ and $b^+$ differ only by a constant.  Clearly, at $q$ the derivative of $\overline{b}^+$ in the direction of $\overline{\gamma}_+'(0)$ is $1$.  Since $\overline{b}^+$ has Lipschitz constant 1, this implies that $\nabla b^+(q) = \overline{\gamma}_+'(0)$.

>From the Bochner formula (\ref{bochner}) and a direct computation one has the following Bochner formula with measure,
\be  \label{Bochner-f}
\frac 12 \Delta_f
|\nabla u|^2 = |\Hess \, u|^2 + \lp \nabla u, \nabla (\Delta_f u) \rp
+ \Ric_f (\nabla u, \nabla u).
\ee
Now apply this to $b^{+}$, since $|| \nabla b^{+}|| = 1$, we have
\be  \label{boch} 0 = || \Hess \, b^{+} ||^2 + \nabla b^{+} ( \Delta_f(b^{+})) +  \Ric_f(\nabla b^{+}, \nabla b^{+}). \ee

Since $\Delta_f(b^{+})$ = 0 and $\Ric_f \geq 0$ we then have that $\Hess \, b^{+} = 0$ which, along with the fact that $||\nabla b^{+}|| =1$ implies  that $M$ splits isometrically in the direction of $\nabla b^{+}$.

To see that $f$ must be constant in the splitting direction note that from (\ref{boch}) we now  have $\Ric_f(\nabla b^{+}, \nabla b^{+}) = 0$ and  $\nabla b^{+}$ points in the splitting direction so  $\Ric(\nabla b^{+}, \nabla b^{+}) = 0$.  Therefore $\Hess f(\nabla b^{+}, \nabla b^{+}) = 0$. Since  $f$ is bounded  $f$ must be constant in $\nabla b^{+}$ direction.
\qed

As Lichnerowicz  points out, the clever covering arguments in \cite{Cheeger-Gromoll1971} along with Theorem \ref{splitting} imply the following structure theorem for compact manifolds with $\Ric_f \geq 0$.

\begin{theo} \label{CG-struc} If $M$ is compact  and $Ric_f \geq 0$ then $M$ is finitely covered by $N \times T^k$ where $N$ is a compact simply connected manifold and $f$ is constant on the flat torus $T^k$.
\end{theo}
Theorem~\ref{CG-struc}  has the following topological consequences.
\begin{coro} \label{Compact-Cor} Let  $M$ be compact with $Ric_f \geq 0$ then
\begin{enumerate}
\item $b_1(M) \leq n$.
\item $\pi_1(M)$ has a free abelian subgroup of finite index of rank $\leq n.$
\item $b_1(M) $ or $\pi_1(M)$ has a free abelian subgroup of rank $n$ if and only if $M$ is a flat torus and $f$ is a constant function.
\item $\pi_1(M)$ is finite if $\Ric_f > 0$ at one point.
\end{enumerate}
\end{coro}

For noncompact manifolds with positive Ricci curvature  the splitting theorem has also been  used by Cheeger and Gromoll \cite{Cheeger-Gromoll1971}  and Sormani \cite{Sormani2001} to give some other topological obstructions.  These results also hold for $\Ric_f$ with $f$ bounded.
\begin{theo} Suppose $M$ is a complete manifold with $\Ric_f > 0$ for some bounded $f$ then
\begin{enumerate}
\item $M$ has only one end
and \item $M$ has the loops to infinity property.
\end{enumerate}
In particular, if $M$ is simply connected at infinity then $M$ is simply connected.
\end{theo}

\sect{Other applications of the mean curvature comparison.} \label{Mean-appl}

Theorem \ref{mean-comp} can be used to prove an excess estimate.  Recall that for $p, q \in M$ the excess function is $e_{p,q}(x) = d(p,x) + d(q,x) - d(p,q)$.  Let $h(x) = d(x, \gamma)$ where $\gamma$ is a fixed minimal geodesic from $p$ to $q$, then (\ref{mean-H=0}) along with the arguments in  \cite[Proposition 2.3]{ Abresch-Gromoll1990}  imply  the following  version of the Abresch-Gromoll excess estimate.

 \begin{theo} \label{excess} Let $\Ric_f \geq 0$,  $|f| \leq k$ and  $h(x) < \min \{d(p,x), d(q,x)\}$ then \[ e_{p,q}(x) \leq 2 \left(\frac{n+4k-1}{n+4k-2}\right) \left(\frac{1}{2} C h^{ n+4k}\right)^{\frac{1}{n+4k-1}}\] where \[ C= 2 \left(\frac{n+4k-1}{n+4k} \right)\left( \frac{1}{d(p,x) - h(x)}+ \frac{1}{d(q,x) - h(x)} \right)  \]
\end{theo}

\Rk (\ref{mean-a}) also implies an excess estimate for manifolds with $\Ric \geq (n-1)H$ and $|\nabla f| \leq a$, however the constant $C$ will depend on $H \cdot d(p,q)^2$ and $e^{ah}$. The mean curvature comparison for $\Ric_f^N$  discussed in the appendix also implies an excess estimate. 

Theorem \ref{excess}  gives extensions of theorems of Abresch-Gromoll \cite{Abresch-Gromoll1990} and Sormani \cite{Sormani2000} to the case where $\Ric_f \geq 0$ for a bounded $f$.

\begin{theo} Let be $M$ a complete noncompact manifold with $Ric_f \geq 0$ for some bounded $f$.
\begin{enumerate}
\item \label {ag} If $M$ has bounded diameter growth and sectional curvature bounded below then $M$ is homeomorphic to the interior of a compact manifold with boundary.
\item \label{sor} If $M$ has sublinear diameter growth then $M$ has finitely generated fundamental group.
\end{enumerate}
\end{theo}

\Rk If we consider $|f|\leq k$, the arguments in \cite{Abresch-Gromoll1990} and \cite{Sormani2000} say slightly more.   Namely,  the diameter growth in  the first part can be of order $\leq \frac{1}{n+4k-1}$ and  in the second part one can derive a explicity constant $S_{n,k}$ such that    the diameter growth only needs to be $\leq  S_{n,k}\cdot r.$  Also see \cite{Yang-Xu-Wang2003} and \cite{Wylie2006}.

We can also apply the mean curvature comparison to the excess function to prove the Myers' theorem.  We note that the excess function was also used to prove a  Myers' theorem in \cite{Petersen-Sprouse1998}.  It is interesting that this proof is exactly suited to our situation, since we only have a uniform bound on mean curvature when $r \leq \frac{\pi}{2 \sqrt{H}}$, while other arguments do not seem to easily generalize.  

\noindent \textbf{Proof of Theorem~\ref{Myers'}} Let $p_1,p_2$ are
two points in $M$ with $d(p_1, p_2) \geq \frac{\pi}{ \sqrt{H}}$ and
set $B = d(p_1,p_2) - \frac{\pi}{\sqrt{H}}$. Let  $r_1(x) =
d(p_1,x)$ and $r_2(x) =d(p_2,x)$ and $e$ be the excess function for the points $p_1$ and $p_2$. By the triangle inequality, $e(x)
\geq 0$ and $e(\gamma(t)) = 0$ where $\gamma$ is a minimal geodesic
from $p_1$ to $p_2$. Therefore, $\Delta_f(e) (\gamma(t)) \geq 0$.

Let $ y_1 = \gamma\left(\frac{\pi}{2\sqrt{H}}\right)$ and $y_2 =
\gamma\left(\frac{\pi}{2\sqrt{H}}+B\right)$. For $i = 1$ and $2$,
$r_i(y_i)= \frac{\pi}{2 \sqrt{H}}$ so, by (\ref{m-f-H>0}), we have
\be \label{Eq-d-1} \Delta_f (r_i(y_i))\leq 2k \sqrt{H}. \ee
(\ref{mean-H<0}) does not give an estimate for  $\Delta_f( r_1
(y_2))$ since $r_1 (y_2)> \frac{\pi}{2 \sqrt{H}}$ but by
(\ref{m_H>0}) and (\ref{Eq-d-1}) we have \be \Delta_f( r_1 (y_2))
\le 2k \sqrt{H} -B(n-1)H.\ee So \be 0 \leq \Delta_f(e)(y_2) =
\Delta_f(r_1)(y_2) + \Delta_f(r_2)(y_2) \le 4k \sqrt{H} -B(n-1)H \ee
which implies $B \le \frac{4k}{(n-1) \sqrt{H}}$ and $d(p_1, p_2) \le
\frac{\pi}{ \sqrt{H}} + \frac{4k}{(n-1) \sqrt{H}}
.$\qed


As we have seen, there is no bound on the distance between two points in a complete manifold with $\Ric_f \geq (n-1)H > 0$.  However, by slightly modifying the argument above one can prove a distance bound between two hypersurfaces that depends on the $f$-mean curvature of the hypersurfaces, here for a hypersurface $N$ the $f$-mean curvature at a point $x \in N$ with respect to the normal vector $n$ is  \be
H^f_{n}(x)= H_{n}(x) +  \langle n, \nabla f \rangle(x) \ee
where $H_{n}$ is the regular mean curvature.   $m_f$ is then the $f$-mean curvature of the geodesic sphere with respect to the inward pointing normal.

\begin{theo}Let $\Ric_f \geq (n-1)H>0$ and let $N_1$ and $N_2$ be two compact hypersurfaces in $M$ then
\be
d(N_1,N_2) \leq \frac{ \max_{p \in N_1} |H_{n_1}^f(p)| + \max_{q \in N_2} |H_{n_2}^f(q)|}{2(n-1)H} \ee
\end{theo}
\Pf
Let $e_{N_1,N_2}(x) = r_1(x)+r_2(x)-d(N_1,N_2)$  where $r_i(x) = d(x, N_i)$.  Then, by applying the Bochner formula to $r_i$ in the same way we applied it to the distance to a point in Section 2, we have
\[ \Delta_f(r_i)(x) \leq  \max_{p \in N_i} |H_{n_i}^f(x)| - (n-1)H d(N_i,x)\]
One now can prove the theorem using a similar argument as in the
proof of Theorem \ref{Myers'}. \qed

A similar argument  also  shows Frankel's Theorem is true for $\Ric_f$.
\begin{theo} Any two compact $f$-minimal hypersurfaces in a manifold with $\Ric_f > 0$ intersect.
\end{theo}
One also has a rigidity statement when $\Ric_f \geq 0$ and $M$ has two $f$-minimal hypersurface which do not intersect, see \cite{Petersen-Wilhelm2003} for the statement and proof in the $\Ric\geq0$ case.

\sect{Examples and Remarks} \label{examples}

The most well known example is the following  soliton, often referred to as the Gaussian soliton.

\begin{Ex} \label{R^n}
Let $M= \mathbb R^n$ with Euclidean metric $g_0$, $f (x)  = \frac{\lambda}{2} |x|^2$. Then $\Hess f =  \lambda g_0$ and $\Ric_f  =   \lambda g_0$.
\end{Ex}
This example shows that, unlike the case of Ricci curvature
uniformly bounded from below by a positive constant, the manifold
could be noncompact  when $\Ric_f \ge \lambda g$ and $\lambda > 0$.

>From this we construct the following.
\begin{Ex} \label{example-H}
Let $M= \mathbb H^n$ be the hyperbolic space. Fixed any $p\in M$,
let $f(x) =(n-1) r^2 = (n-1)d^2(p,x)$. Now Hess $r^2=  2|\nabla r|^2 + 2 r \Hess r \ge 2I $, therefore
$\Ric_f \ge (n-1)$.
\end{Ex}
This example shows that the Cheeger-Gromoll splitting theorem and Abresch-Gromoll's excess estimate do not hold for $\Ric_f \ge 0$, in fact they don't even hold for $\Ric_f \geq \lambda > 0$.  Note that the only properties of hyperbolic space used are that $\Ric \geq -(n-1)$ and that Hess $r^2 \geq 2I.$  But $\Hess \, r^2 \geq 2I$ for any Cartan-Hadamard manifold, therefore any Cartan-Hadamard manifold with Ricci curvature bounded below has a metric with $\Ric_f \geq 0$ on it.   On the other hand, in these examples $\Ric < 0$. When $\Ric < 0$ ($\Ric \le 0$) and $\Ric_f \ge 0 (\Ric_f >0)$, then $\Hess f  >0$ and $f$ is strictly convex. Therefore $M$ has to homeomorphic to $\mathbb R^n$.

A large class of examples are given by gradient Ricci solitons.
Compact expanding or steady solitons are Einstein ($f$ is constant)
\cite{Perelman-math.DG/0211159}. There are nontrivial compact shrinking solitons
\cite{Koiso1990,Cao-Huai-Dong1996}. These examples also have positive Ricci curvature but in the noncompact case there are Kahler Ricci shrinking solitons that do not have nonnegative Ricci curvature \cite{Feldman-Ilmanen-Knopf2003}.
 Clearly
there are examples with $\Ric_f \ge 0$ but Ricci curvature is not
nonnegative, like Example~\ref{example-H}. One can also construct
example that $f$ is bounded. In fact one can use the following
general local perturbation.

\begin{Ex}  \label{local-perturbation}
Let $M^n$ be a complete Riemannian manifold with $\Ric \ge 0$ except in a neighborhood of a point $p, U_p$. If the sectional curvature is $\le 1$ and $\Ric \ge -1$ on $U_p$ and $U_p \subset B(p, \pi/4)$, then the distance function $r (x) = d(p,x)$ is strictly convex (with a singularity at $p$) on $U_p$. On $B(p,\pi/4)$ let $f =  r^2$. If $\Ric \ge 1$ on the annulus $A(p,\pi/4,\pi/2)$, then one can extend $f$ smoothly to $M$ so that $f$ is constant outside $B(p, \pi/2)$ and $\Ric_f \ge 0$ on $M$.
\end{Ex}

The following example shows that there are manifolds with $Ric_f \geq 0$ which do not have polynomial $f$-volume growth.

\begin{Ex} \label{example-R}
Let $M= \mathbb R^n$ with Euclidean metric, $f(x_1, \cdots, x_n) =
x_1$. Since Hess $f= 0$, $\Ric_f = \Ric = 0$. On the other hand $\vol_f (B(0,r))$ is of exponential growth.  Along the $x_1$
direction, $m_f - m_H = -1$ which does not goes to zero.
\end{Ex}

In this example $|\nabla f| \le 1$, so $\Ric_f \ge 0$ and $|\nabla
f|$ bounded does not imply polynomial  $f$-volume growth either.

\begin{ques}  \label{Ric-Ric_f}
If $M^n$ has a complete metric  and measure such that $\Ric_f \ge 0$ and $f$ is bounded, does $M^n$ has a metric with $\Ric \ge 0$?
\end{ques}

There is no counterexample even without the $f$ bounded condition.

It is also natural to consider the scalar curvature with measure. As pointed out by Perelman in \cite[1.3]{Perelman-math.DG/0211159} the corresponding scalar curvature equation is $S_f = 2 \Delta f - |\nabla f|^2 + S$.  Note that this is different than taking the  trace of $Ric_f$ which is $\Delta f + S$. However, The Lichnerovicz formula and theorem naturally extend to $S_f$. But $\Ric_f \ge 0$ doesn't imply $S_f \ge 0$ anymore.  So one can ask the following question.
\begin{ques}
If $M^n$ is a compact spin manifold with $\Ric_f >0$, is the $\hat{A}$-genus zero?
\end{ques}

One could  try to see if the $K3$ surface has a metric with $\Ric_f >0$. If this were true it would  give a negative answer to Question~\ref{Ric-Ric_f}.

\appendix

\section{Mean curvature comparison for
$N$-Bakry-Emery Ricci tensor} \label{N-dim-appendix}
\renewcommand{\theequation}{A-\arabic{equation}}

In \cite{Bakry-Qian2005} the volume comparison theorem and Myers'
theorem  for the $N$-Bakry-Emery Ricci tensor are proven using what
we have called  a mean curvature comparison (actually their work is
slightly more general than the cases treated in this paper).     In
this appendix, for clarity, we state this comparison in the language
we have used above, which is standard in Riemannian geometry.

Recall the definition of the $N$-Bakry-Emery tensor is \[ \Ric_f^N =
\Ric_f -\frac 1N df \otimes df  \ \ \  \mbox{for} \  N>0.
 \]

 The main idea is that the a Bochner formula holds for $\Ric_f^N$ that looks like the Bochner
 formula for the Ricci tensor of an $n+N$ dimensional manifold .  This formula seems to have been
 Bakry and Emery's original motivation for the definition of the Bakry-Emery Ricci tensor and for
 their more general curvature dimension inequalities for diffusion operators \cite{Bakry-Emery1985}.
 See \cite{Ledoux2000, Li2005} for elementary proofs of the inequality.

\[\frac 12 \Delta_f
|\nabla u|^2 \geq \frac{(\Delta_f(u))^2}{N+n} + \lp \nabla u, \nabla (\Delta_f u) \rp
+ \Ric_f^N (\nabla u, \nabla u) \]

For the distance function, we actually have \[ m'_f \leq -\frac{(m_f)^2}{n+N-1} - \Ric^N_f(\partial_r, \partial_r).\]
  Thus, using the standard Sturm-Liouville comparison argument, or an argument similar to the one we give above,  one has the  mean curvature comparison.
\begin{theo} [Mean curvature comparison for $N$-Bakry-Emery] \cite{Bakry-Qian2005}
Suppose that $N > 0$ and $\Ric^N_f \geq (n+N-1)H$, then \[ m_f(r) \leq m_H^{n+N}(r).\]
\end{theo}

This comparison along with the methods used above gives proofs of the comparison theorems for $\Ric_f^N$.

The Bochner formula  for $\Ric_f^N$ has many  other applications to other geometric problems not treated here such as eigenvalue problems and Liouville theorems, see for example \cite{Bakry-Qian2000} and \cite{Li2005} and the references there in.

 In \cite{Lott2003}  Lott  shows that if $M$ is compact with $\Ric_f^q \geq \lambda$ for some positive integer $q \ge 2$,
 then, in fact,  there is a family of  metrics on $M \times S^q$ with Ricci curvature bounded below by $\lambda$. Moreover, the metrics on the sphere collapse so that $M$ is a Gromov-Hausdorff limit of $n+q$ dimensional manifolds with  Ricci curvature bounded below by $\lambda$.  This  gives an alternate approach to prove many of the comparison  and topological theorems for $\Ric_f^q$.

\setcounter{equation}{0}
\section{ODE proof of mean curvature comparison}
 \renewcommand{\theequation}{B-\arabic{equation}}

\begin{theo}[Mean Curvature Comparison] \label{Old-mean-comp} Assume $\Ric_f (\partial_r, \partial_r) \ge (n-1)H$,

a) if $\partial_r f \ge -a \ (a \ge 0)$ along a minimal geodesic
segment from $p$ (when $H>0$ assume $r \le \pi/2\sqrt{H}$) then
 \be \label{mH<0}
m_f(r) - m_H(r) \le a
 \ee
 along that minimal geodesic segment from $p$. Equality holds if and only if the radial sectional
 curvatures are equal to $H$ and $f(t) = f(p) - at$ for all $t<r$.
\\In particular when $a= 0$, we have
 \be
 m_f(r) \le m_H(r)
 \ee
 and equality holds if and only if all radial sectional curvatures are $H$ and $f$ is constant.

 b) if $|f| \le k$ along a minimal geodesic
segment from $p$ (when $H>0$ assume $r \le \pi/2\sqrt{H}$) then
 \be  \label{old-mean-H<0}
m_f(r) - m_H  \le (n-1) e^{\frac{4k}{n-1}} \left(\frac{\sqrt{|H|} \sn_H (2r) + 2|H|r}{\sn_H^2(r)} \right)
  \ee
  along that minimal geodesic segment from $p$, where $\sn_H(r)$ is the unique function satisfying \[
  \sn_H ''(r) +H \sn_H(r) = 0, \ \ \sn_H(0) = 0, \ \ \sn_H'(0) = 1.\]
In particular when $H=0$ we have
\be \label{old-mean-H=0} m_f(r) - \frac{n-1}{r}  \le 4(n-1)
e^{\frac{4k}{n-1}}\frac 1r. \ee
\end{theo}

\Pf We compare $m_f$ to the mean curvature of the model space. Note
that the mean curvature of the model space $m_H$ satisfies \be
\label{Riccati-inequ-H} m'_H = -\frac{m_H^2}{n-1} - (n-1)H. \ee
Using $\Ric_f \ge (n-1)H$, and subtracting (\ref{Riccati-inequ-mf})
by (\ref{Riccati-inequ-H}) gives \ba
 (m_f-m_H)'  & \leq & -\frac{1}{n-1} \left[ (m_f+\partial_r f)^2 - m_H^2 \right]  \label{premc} \\
& = & -\frac{1}{n-1} \left[ (m_f-m_H+\partial_r f)(m_f+m_H +\partial_r f) \right]. \label{mc}
\ea

\noindent {\bf Proof of Part a)}:  Write (\ref{mc}) as the following \be
 (m_f-m_H -a)' \le -\frac{1}{n-1} \left[ (m_f-m_H-a + a +\partial_r
f)(m_f+m_H +\partial_r f) \right]. \label{mc-b} \ee Let us define
$\psi_{a,H} = \max\{m_f-m_H -a,0\} = (m_f-m_H-a)_+$, and declare
that $\psi_{a,H}=0$  whenever it becomes undefined. Since
$\partial_r f \ge -a$, $a+\partial_r f \ge 0$. When $\psi_{a,H}
\ge 0,\ m_f+m_H +\partial_r f \ge a+ \partial_r f +2m_H  \ge 2m_H$
which is $ \ge 0$ if $m_H \ge 0$. Using this and  (\ref{mc-b}) gives
 \be
\label{ode} \psi'_{a,H} \le -\frac{1}{n-1}  (m_f+m_H +\partial_r f)
\psi_{a,H} \le 0. \ee
Since $\lim_{r \ra 0} \psi_{a,H} (r) = (-\partial
f_r (0) -a)_+ = 0$, we have
$\psi_{a,H} (r) = 0$ for all $r \ge 0$. This finishes the
proof of the inequality.

Now suppose that $m_f = m_H + a$, then from (\ref{premc}) we have that $m = m_H$ which implies that $\partial_r f = -a$. So $\partial_r^2f \equiv 0$ which then implies that $\Ric(\partial_r, \partial_r)= \Ric_f(\partial_r, \partial_r) \geq (n-1)H$. Now the rigidity for the regular mean curvature comparison implies that all the sectional curvatures are constant and equal to $H$.
\newline


\noindent {\bf Proof of Part b)}: Write (\ref{mc}) as \ba
(m_f-m_H)' & \le & -\frac{1}{n-1} \left[ (m_f-m_H+\partial_r f)(m_f - m_H +2m_H +\partial_r f) \right]\nonumber \\
& = & -\frac{1}{n-1} \left[ (m_f-m_H)^2 +2(m_H+\partial_r
f)(m_f-m_H) + \partial_r f (2m_H+\partial_r f)\right]. \label{mmh}\ea
Now define $\psi = \max\{m_f-m_H,0\} = (m_f-m_H)_+$, the error from the mean curvature
 comparison, and declare that $\psi = 0$ whenever it becomes
 undefined. Define
  \be \label{rho}
\rho = \left[ -\frac{1}{n-1}\partial_r f (2m_H+\partial_r f)\right]_+.
\ee
Using this notation and inequality (\ref{mmh}) we obtain
 \be \label{main-ineq}
 \psi' \leq -\frac{1}{n-1} \psi^2 -\frac{2}{n-1}(m_H+\partial_rf)\psi +\rho.
 \ee
We would like to estimate $\psi$ in terms of $\rho$. It is enough to
use the linear differential inequality
\be \label{linear}
 \psi' + \frac{2}{n-1}(m_H+\partial_rf)\psi \le \rho.
 \ee
When $\partial_r f = 0$ (in the usual case), we have $\rho = 0$ and $\psi = 0$, getting the classical mean curvature comparison. In general,
by (\ref{rho}), the
definition of $\rho$, when $m_H>0$
  \be \label{rho>0}
  \rho > 0 \Longleftrightarrow  -2m_H < \partial_r f < 0.
  \ee
Also
\be  \label{rho<}
\rho \le \left( -\frac{2} {n-1}(\partial_r f) m_H\right)_+ .
\ee
 Therefore we have
\[
\rho  \le \frac{4}{n-1} m_H^2.
  \]
Note that  $m_H = (n-1) \frac{\sn_H'(r)}{\sn_H (r)}$. Now (\ref{linear}) becomes
\[ \psi' +\left(2  \frac{\sn_H'(r)}{\sn_H (r)}  +
\frac{2}{n-1}\partial_rf\right)\psi \le  4 (n-1)  \left(\frac{\sn_H'(r)}{\sn_H (r)}\right)^2.
\]
 Multiply  this
 by the integrating factor  $\sn_H^2(r)e^{\frac{2}{n-1} f(r)}$ to obtain
\[ \left( \sn_H^2(r) e^{\frac{2}{n-1} f(r)} \psi(r)\right)' \le 4(n-1) e^{\frac{2}{n-1} f(r)} (\sn_H'(r))^2.\]
Since $\psi (0)$ is bounded, integrate this from $0$ to $r$ gives
\be  \label{sn}
\sn_H^2(r) e^{\frac{2}{n-1} f(r)} \psi(r) \le 4(n-1) \int_0^r
e^{\frac{2}{n-1} f(t)} (\sn_H'(r))^2 dt.\ee
Since $|f| \le k$, we have
\[
\psi (r) \le (n-1) e^{\frac{4k}{n-1}} \left(\frac{\sqrt{|H|} \sn_H(2r) + 2|H|r}{ \sn_H^2(r)} \right).
  \]
When $H = 0$, $\sn_H(r) = r$, from (\ref{sn})
we get
  \[\psi (r) \le 4(n-1) e^{\frac{4k}{n-1}}\frac{1}{r}.
  \]
  This completes the proof of Part b).
  \qed


Department of Mathematics,

University of California,

Santa Barbara, CA 93106

wei@math.ucsb.edu
\newline

Department of Mathematics,

University of California,

Los Angeles, CA

wylie@math.ucla.edu
\end{document}